\input amstex
\documentstyle{amsppt}
\document
\topmatter
\title
Bi-Hermitian Gray surfaces II.
\endtitle
\author
Wlodzimierz Jelonek
\endauthor

\abstract{The aim of this paper is to classify bi-Hermitian
compact surfaces $(M,g)$ whose Ricci tensor $\rho$  satisfies the
relation $\nabla_X\rho(X,X) =\frac13X\tau g(X,X)$.  }
 \endabstract

\endtopmatter
\define\G{\Gamma}
\define\DE{\Cal D^{\perp}}
\define\e{\epsilon}
\define\n{\nabla}
\define\om{\omega}
\define\w{\wedge}
\define\k{\diamondsuit}
\define\th{\theta}

\define\a{\alpha}
\define\be{\beta}

\define\lb{\lambda}
\define\A{\Cal A}
\define\AC{\Cal AC^{\perp}}
\define\bO{\bar \Omega}
\define\1{\Cal D_{\lb}}
\define\2{\Cal D_{\mu}}
\define\0{\Omega}

\define\bt{\bar\tau}

\define\bth{\bar{\th}}
\define\bJ{\bar J}
\define\bg{\bar g}
\define\bn{\overline \nabla}
\define\brho{\bar\rho}

\define\Si{\Sigma}
\define\J{\Cal J}
\define\De{\Cal D}
\define\de{\Delta}

\define\dl{\delta}

\define\m{(M,g,J)}
\define\pe{(M,g)}
\centerline{\it This paper I  dedicate to the memory of Alfred Gray.}

\bigskip
{\bf 0. Introduction. } Let $(M,g)$ be a Riemannian manifold with the Ricci tensor $\rho$ which   satisfies the condition
$$\n_X\rho(X,X) =\frac2{n+2}X\tau g(X,X)\tag *$$
where $\tau$ is the scalar curvature of $(M,g)$ and $n=$dim$M$. There are many interesting manifolds which satisfy (*). Among them are (compact)
 Einstein-Weyl manifolds, weakly self-dual K\"ahler surfaces (see [J-1],[J-2] and [A-C-G])  and D'Atri spaces.
The  property (*) was studied by A. Gray in [G] (see also [Be] p.433). A. Gray called Riemannian manifolds satisfying $(*)$ the $\AC$ manifolds.
In [J-1] we  showed that every K\"ahler surface has a harmonic anti-self- dual part $W^-$ of the Weyl tensor $W$ (i.e. such that $\dl W^-=0$) if and
only if it is an  $\AC$-manifold.   In [J-1] we have also showed that any simply connected 4-dimensional $\AC$-manifold $(P,g)$ , whose Ricci tensor
has exactly two eigenvalues of multiplicity 2, admits two opposite to each other Hermitian structures which commute with the Ricci tensor.

It is not difficult to prove that a compact 4-manifold with even first Betti number
admitting two opposite to each other Hermitian structures $J,\bJ$ which commute with the Ricci tensor $\rho$ of $(P,g)$ is a ruled surface or is
locally a product of two Riemannian surfaces [see [J-2]).
 In [J-2] we have given the example of a K\"ahler $\AC$-metric on a Hirzebruch surface $F_1$ (which was also independently constructed in [A-C-G])
  and in [J-3] we have constructed families of bi-Hermitian Gray surfaces on all the Hirzebruch surfaces $F_k$. These are, apart from two exceptional
   families of metrics on $F_1$ and one exceptional family on $F_2$,  all co-homogeneity one bi-Hermitian Gray metrics on ruled surfaces of genus $g=0$.

The aim of the present paper is to describe compact $\AC$-4-manifolds $(M,g)$ with non-constant scalar curvature, admitting two oppositely definite
 Hermitian structures $J,\bJ$ commuting with the  Ricci tensor of $(M,g)$.  We shall call such surfaces the bi-Hermitian Gray surfaces. Surfaces which
  admit two
oppositely oriented complex structures will be called the bi-Hermitian surfaces. We should warn the reader that the notion of a bi-Hermitian surface
has been recently used also in the different context (see [A-G-G] where a bi-Hermitian surface means a surface admitting two positively oriented Hermitian
 structures).
We show in the present paper that compact  bi-Hermitian Gray surfaces  with non constant scalar curvature and even first Betti number  are ruled surfaces
which  ( at least if their genus $g\ge 1$) are local cohomogeneity 1 with respect to the local group of all local isometries of $(P,g)$. We shall
classify all bi-Hermitian Gray surfaces which are ruled surfaces of genus $g>0$. At first we shall show that there exists an open and dense subset
$U$ of $P$ such that $U=(a,b)\times P_0$ where $P_0$ is a three dimensional $\A$-manifold which is a circle bundle over a compact Riemann surface of
 constant sectional curvature. Consequently if $P_0$ is not a trivial bundle then it coincides with the space $G\slash\G$ where $G$ is one of the
  groups:$ SU(2),H,\widetilde{ SL(2,\Bbb R)}$ where $H$ means the Heisenberg group and $\G$ is a discrete subgroup of  $Iso(G)$.
Using the methods of B. Bergery (see [B],[S]) we reduce the problem to a certain ODE of the second order. We shall find all  positive solutions of
 these equations satisfying the appropriate boundary conditions.  In this way we classify compact bi-Hermitian Gray surfaces of genus $g>0$ and also
 give  new examples of  compact 4-dimensional $\AC$-manifolds (problem of finding such manifolds was stated in  [Be] p.433).  In the last section of
 the paper we
describe in an explicit way co-homogeneity one $\AC$ - metrics on $\Bbb{CP}^2 $ whose Ricci tensor  is invariant with respect to the standard complex
structure $J$ of $\Bbb{CP}^2$ and such that the opposite Hermitian structure $\bJ$ is defined on $\Bbb{CP}^2-\{x_0\}$ for some $x_0\in\Bbb{CP}^2$.

\par
\bigskip
 {\bf 1. Hermitian 4-manifolds. }     Let $(M,g,J)$ be an almost Hermitian manifold, i.e. $(M,g)$ is a Riemannian manifold
 and $J:TM\rightarrow TM$ satisfies $J^2=-id_{TM}$ and $g(JX,JY)=g(X,Y)$ for all $X,Y\in TM$. We say that $(M,g,J)$ is a Hermitian  manifold if its almost Hermitian
  structure $J$ is integrable, i.e. $J$ is an orthogonal complex structure. In the sequel we shall consider 4-dimensional Hermitian manifolds $(M,g,J)$
  which we shall also call Hermitian surfaces.
  Such manifolds are always oriented and we choose an orientation in such a way that the K\"ahler form $\0(X,Y)=g(JX,Y) $ is a self-dual form
  (i.e. $\0\in\w^+M$). The vector bundle of self-dual forms admits a decomposition
$$\w^+M=\Bbb{R}\0\oplus LM,\tag 1.1$$
where by $LM$ we denote the bundle of real $J$-skew invariant 2-forms (i.e $LM=\{\Phi\in\w M:\Phi(JX,JY)=-\Phi(X,Y)$\}). The bundle $LM$ is a complex
line bundle over $M$ with the complex structure $\Cal J$ defined by $(\Cal J\Phi)(X,Y)=-\Phi(JX,Y)$.
For a 4-dimensional Hermitian manifold the covariant derivative of the K\"ahler form $\0$ is locally expressed by
$$\n \0=a\otimes\Phi+\J a\otimes\J\Phi, \tag 1.2$$
where $\J a(X)=-a(JX)$.
The Lee form $\th$ of $(M,g,J)$ is defined by the equality
$$d\0=\th\w \0\tag 1.3$$
 We have $\th=-\delta\0\circ J$.
 A Hermitian manifold $(M,g,J)$ is said to have Hermitian Ricci tensor $\rho$ if $\rho(X,Y)=\rho(JX,JY)$ for all $X,Y\in \frak X(M)$.
An opposite (almost)  Hermitian structure  on a Hermitian 4-manifold $(M,g,J)$ is an (almost) Hermitian structure $\overline J$ whose  K\"ahler form
( with respect to $g$) is anti-self-dual.

A distribution $\De\subset TM$ is called umbilical  if $\n_XX_{|\DE}=g(X,X)\xi$ for every $X\in\G(\De)$, where $X_{|\DE}$ is the $\DE$ component of $X$ with
 respect to the orthogonal decomposition $TM=\De\oplus\DE$. The vector field $\xi$ is called the mean curvature normal of $\De$.  An involutive
 distribution $\De$ is tangent to a  foliation, which  is called totally geodesic if its every leaf is a totally geodesic submanifold of $(M,g)$
  i.e. $\n_XX\in \De$ if $X$ is a section of a vector bundle $\De\subset TM$. In the sequel we shall not distinguish between $\De$ and a tangent
   foliation and we shall also say that $\De$ is totally geodesic in such a case.

On any  Hermitian non-K\"ahler 4-manifold $\m$ there are  two natural distributions $\De=\{X\in TM:\n_XJ=0\}$,   $\DE$ defined in the open set
$U=\{x:|\n J_x|\ne 0\}$. The distribution $\De$ we shall call the nullity distribution of $(M,g,J)$. From  (1.2) it is clear that $\De$ is
$J$-invariant and that dim$\De=2$ in $U=\{x\in M:\n J_x\ne 0\}$. By $\DE$ we shall denote the orthogonal complement of $\De$ in $U$.  On $U$
we can define the opposite almost Hermitian structure $\bJ$ by formulas $\bJ X=JX$ if $X\in \DE$ and $\bJ X=-JX$ if $X\in\De$ which we shall
call natural opposite almost Hermitian structure. It is not difficult to check that for the famous Einstein Hermitian manifold
 $\Bbb{CP}^2\sharp\overline{\Bbb{CP}}^2$ with Page metric (see [P],[B],[S],[K],[LeB]) the opposite structure $\bJ$ is Hermitian and this structure
  extends  to the global opposite Hermitian structure.

  A ruled surface of genus $g$ is a complex surface $X$ admitting
  a ruling, i.e. an analytically  locally trivial fibration with
  fibre $\Bbb CP^1$ and structural group $PGL(2,\Bbb C)$ over a
  smooth compact complex curve (a Riemannian surface) of genus
  $g$.

\par
\bigskip

By an $\AC$- manifold (see [G]) we mean a Riemannian manifold $(M,g)$ satisfying the condition
$$\frak C_{X Y Z}\n_X\rho(Y,Z)=\frac 2{(\text{dim}M+2)}\frak C_{X Y Z}X\tau g(Y,Z),  \tag 1.4$$
where $\rho$ is the Ricci tensor of $(M,g)$ and $\frak C$ means the cyclic sum. A Riemannian manifold $(M,g)$ is an $\AC$ manifold if and only if the
 Ricci endomorphism $Ric$ of $(M,g)$ is of the form $Ric=S+\frac{2}{n+2}\tau Id$ where $S$ is a Killing tensor, $\tau$ is the scalar curvature and
 $n=$dim$M$. Let us recall that a (1,1) tensor $S$ on a Riemannian manifold $(M,g)$ is called a Killing tensor if $g(\n S(X,X),X)=0$ for all $X\in TM$.
 Let us recall a result from [J-1]:
 \bigskip
{\bf Lemma 0. } {\it Let  $S$ be a Killing tensor on a
4-dimensional Riemannian manifold $(M,g)$. Let us assume that $S$
has two 2-dimensional oriented eigendistributions $D_1,D_2$. Then
there exist two opposite Hermitian complex  structures $J,\bJ$ on
 $M$ which commute with $S$. }
 \bigskip
 It is not difficult to prove the following lemmas:
\medskip
{\bf Lemma  1. } {\it Let $S\in End(TM)$ be a (1,1) tensor on a
Riemannian 4-manifold $(M,g)$. Let us assume that $S$ has exactly
two everywhere different eigenvalues $\lb,\mu$ of the same
multiplicity 2,  i.e. \text{ dim} $\1$=\text{ dim }$\2=2$, where
$\1,\2$ are eigendistributions of $S$ corresponding to $\lb,\mu$
respectively. Then $S$ is a Killing tensor if and only if both
distributions $\1$ and $\2$ are umbilical with mean curvature
normal equal respectively}
$$\xi_{\lb}=\frac{\n\mu}{2(\lb-\mu)}, \  \xi_{\mu}=\frac{\n\lb}{2(\mu-\lb)}.$$

\medskip
{\bf Lemma 2. } {\it Let $(M,g)$ be a 4-dimensional Riemannian
manifolds whose Ricci tensor $\rho$ has two eigenvalues
$\lb(x),\mu(x)$ of the same multiplicity 2 at every point $x$ of
$M$. Let us assume that the eigendistribution $\De_{\lb}=\De$
corresponding to $\lb$ is a totally geodesic foliation and the
eigendistribution $\De_{\mu}=\De^{\perp}$ corresponding to $\mu$
is umbilical. Then  $(M,g)$ is an $\AC$-manifold if and only if
$\lb-2\mu$ is constant and $\n\tau\in\G(\De)$. The distributions
$\De,\DE$ determine two Hermitian structures $J,\bJ$ which are
opposite to each other and commute with $\rho$. Both structures
$J,\bJ$ are Hermitian complex  and $\De$ is contained in the
nullity of $J$ and $\bJ$.}
 \bigskip
In the sequel we shall need the following two lemmas.
\medskip
{\bf Lemma   A. } {\it Let assume that $(M,g,J)$ is a compact Hermitian c.K. surface with
Hermitian Ricci tensor $\rho$. If $\zeta$ is a holomorphic Killing vector field on $(M,g,J)$,
then $\th(\zeta)=0$, where $\th$ is a Lee form of $(M,g,J)$.}
\medskip
{\it Proof.} Let $\0$ be a K\"ahler form of $(M,g,J)$. Then $L_{\zeta}\0=0$. Since $d\0=\th\w\0$
and $d\th=0$ it follows that $0=d(L_{\zeta}\0)=L_{\zeta}(d\0)=L_{\zeta}(\th\w\0)=L_{\zeta}(\th)\w\0$. Thus $L_{\zeta}(\th)=0$. Consequently
$d(\th(\zeta))=0$. It follows that $\th(\zeta)$ is constant on $M$, and consequently equals $0$, since the set $\{x\in M:|\th|_x=0\}$ is non-empty
 (see [J-4]).$\k$
\medskip
{\bf Lemma B.} {\it Let us assume that $S$ is a Killing tensor on four dimensional manifold $(M,g)$ with two eigenvalues everywhere distinct and with
 two-dimensional oriented eigen-distributions. Let $J,\bJ$ be Hermitian structures on $(M,g)$ determined by $S$. Let $\xi$ be a Killing vector field on
 $(M,g)$ such that $L_{\xi}S=L_{\xi}J=L_{\xi}\bJ=0$ and $\n_{\xi}J=\n_{\xi}\bJ=0$.  Then $S\xi$ is a Killing vector field on $(M,g)$.}

\medskip
{\it Proof. }  Let us define $TX:=\n_X\xi$.  Then $T\circ J=J\circ T$ and analogously $T\circ \bJ=\bJ\circ T$. Let us define $p=J\circ \bJ$.
It is clear that
$p\circ T= T\circ p.$
Consequently $$S\circ T=T\circ S.\tag 1.5$$
From $(1.5)$ we obtain $\n_{\xi}S=0$. Now we shall show that a field $\zeta=S\xi$ is Killing.
We have
$$\gather
g(\n_X\zeta,X)= g(\n S(X,\xi),X)+g(S(\n_X\xi),X)=\\=-\frac12g(\n_{\xi}S(X),X)
+g(STX,X)=g(STX,X)=0
\endgather$$
since $2g(\n S(X,\xi),X)+g(\n_{\xi}S(X),X)=0$ because $S$ is a Killing tensor and
 $$g(STX,X)=  g(TSX,X)=  -g(TX,SX)=  -g(STX,X).\k$$
We also have (see [J-4])
\medskip
{\bf Theorem 0.} {\it Let us assume that $(M,g,J)$ is a  compact
conformally K\"ahler non-K\"ahler Hermitian surface and let
$(M,g_0,J)$ be a K\"ahler surface in the conformal class
$(M,[g],J)$. Then both $(M,g,J), (M,g_0,J)$ admit a holomorphic
Killing vector field $\xi$ with zeros, such that $\n_{\xi}J=0$,
where $\n$ is the Levi-Civita connection of $(M,g)$. Moreover
$\xi=J\bn u$ where $\bn$ is the Levi-Civita connection of
$(M,g_0)$, $u$ is a positive, smooth function on $M$ such that
$g=u^{-2}g_0$ and we have
$$\n_{\xi_0}\xi_0=-\n\a+\a J\xi_0,\hskip0.5cm  g_0(\xi,\xi)=\a^2,\hskip0.5cm \th(X)=-2\a g(J\xi_0,X)\tag 1.4$$
where $\xi_0=\frac1{\sqrt{g(\xi,\xi)}}\xi$ is defined in
$U=\{x:\xi_x\ne 0\}$  and $\a=\frac1{2\sqrt2}\sqrt{g(\n J,\n J)}$.
Killing vector field $\xi$
 belongs to the center $\frak z(\frak{iso}(M))$ of the Lie algebra $\frak{iso}(M)$ of the group of isometries $Iso(M)$ of $(M,g)$.}
\bigskip

{ \bf 2. Bi-Hermitian Gray surfaces. }  Let us assume that $(M,g)$
is a compact irreducible 4-dimensional $\AC$-manifold whose Ricci
tensor has two eigenvalues $\lb,\mu$. The following definition we
shall use in he sequel.
\medskip
{\it Definition. } Bi-Hermitian Gray  surface is  an irreducible
$\AC$ 4-manifold $(M,g)$,  which admits two different Hermitian
complex  structures $J,\bJ$ of opposite orientation which commute
 with the Ricci tensor $\rho$ of $(M,g)$.
\medskip
Let $\n$  be a Levi-Civita connection and $\rho$ the Ricci tensor of $(M,g)$.
 We say that $(M,g)$ is a proper bi-Hermitian Gray surface if $\n\rho\ne 0$ or equivalently if the scalar curvature $\tau$ of $(M,g)$ is non-constant.
We shall assume in the sequel that $(M,g)$ is a bi-Hermitian Gray
surface with even first Betti number ($b_1(M)$ is even). Due to
the results  of Apostolov and Gauduchon [A-G-1] it  follows that
both Hermitian surfaces $(M,g,J)$ and $(M,g,\bJ)$ are locally
conformally K\"ahler, hence they are both
 conformally K\"ahler due to a result by I. Vaisman [V]. Thus there exist a K\"ahler surface $(M,\bg,J)$   and a K\"ahler surface $(M,\bg_1,\bJ)$
 which are conformally equivalent to $(M,\bg)$ where $\bg_1=h^2g,\bg=f^2g$  and $f,h$ are some smooth functions on $M$. In our paper [J-1] we
 have proved that an oriented  4-dimensional Riemannian manifold $(M,g)$ admitting a Killing tensor $S\in End(TM)$ with exactly two eigenvalues
 $\lb,\mu$ everywhere distinct admits (up to two-fold covering) two Hermitian structures commuting with $S$  and oppositely oriented
 (see Prop.3 in [J-1] and Lemma 0 in the present paper). Hence every 4-dimensional $\AC$-manifold,  whose Ricci tensor has two everywhere distinct eigenvalues admits
 (up to four fold covering) two oppositely oriented Hermitian structures commuting with the Ricci tensor $\rho$ of $(M,g)$. Now we prove

\bigskip
{\bf Proposition 1. } {\it Let us assume that $(M,g)$ is a compact irreducible  bi-Hermitian Gray surface with even first Betti number. Then $(M,g)$ is
an Einstein Hermitian   manifold $\Bbb{CP}^2\sharp\overline{\Bbb{CP}}^2$ with D. Page's metric or the eigenvalues of the Ricci tensor of $(M,g)$ are
everywhere distinct.}
\medskip
{\it Proof.} Let us denote by $J,\bJ$ the opposite Hermitian structures on $(M,g)$ such that $S\circ J=J\circ S$,$S\circ \bJ=\bJ\circ S$ where $S$ is the
 Ricci tensor of $(M,g)$. Let $\{E_1,E_2,E_3,E_4\}$ be a local  orthonormal frame on $(M,g)$ such that $E_1,E_2\in\1$, $E_3,E_4\in\2$ where $\1,\2$
  are eigensubbundles of $S_0$ and
$$JE_1=E_2,\bJ E_1=E_2,JE_3=E_4,\bJ E_3=-E_4.$$
Since $(M,g)$ is an $\AC$-manifold it follows that $S=S_0+\frac{\tau}3Id$ where  $S_0$ is a Killing tensor on $(M,g)$ (we identify $(1,1),(2,0),(0,2)$
tensors on $(M,g)$ by means of $g$). From [J-1] (2.21) it follows that
$$\gather   (\mu-\lb)(\n J(E_1,E_1)+\n J(E_2,E_2))=-J(\n\lb)+(S_0-\lb)([E_1,E_2]),\tag 2.1a\\
(\lb-\mu)(\n J(E_3,E_3)+\n J(E_4,E_4))=-J(\n\mu)+(S_0-\mu)([E_3,E_4]),\tag2.1b\endgather$$
where $\lb,\mu$ are eigenvalues of $S_0$. Consequently
$$(\mu-\lb)(tr_g\n J)=J(\n\mu-\n\lb)+(S_0-\lb Id)([E_1,E_2])-(S_0-\mu Id)([E_3,E_4]).$$

Thus
$$(tr_g\n J)=J(\n\ln|\mu-\lb|)+[E_1,E_2]_{\mu}+[E_3,E_4]_{\lb},\tag 2.2$$
where by $X_{\lb},X_{\mu}$ we mean the components of $X\in TM$ with respect to the decomposition $TM=\1\oplus \2$. Hence in the set
$U=\{x\in M:\lb(x)\ne\mu(x)\}$ the following relation holds
$$ -J(tr_g\n J)=\n\ln|\mu-\lb|-J([E_1,E_2]_{\mu}+[E_3,E_4]_{\lb}).$$
Analogously one can prove that in $U$
$$ -\bJ(tr_g\n \bJ)=\n\ln|\mu-\lb|-\bJ([E_1,E_2]_{\mu}-[E_3,E_4]_{\lb}).$$
The above equations yield that in $U$
$$\th+\bth=2d\ln|\mu-\lb|,\tag 2.3$$
where $\th,\bth$ are the Lee forms of $(M,g,J),(M,g,\bJ)$ respectively. Since $b_1(M)$ is even it follows that both surfaces $(M,g,J),(M,g,\bJ)$ are
 conformally K\"ahler. It means that there exist smooth, positive functions $f,h\in C^{\infty}(M)$ such that $(M,f^2g,J)$,  $(M,h^2g,\bJ)$ are K\"ahler.
 Consequently
$\th=-2d\ln f,\bth=-2d\ln h$. Thus there exists a constant $C\in \Bbb R-\{0\}$ such that
 $$fh=\frac C{\lb-\mu}.\tag 2.4$$
Now $f,h$ are globally defined, smooth functions on $M$ hence they are bounded. It follows that
$U=M$ or $U=\emptyset$. Since $(M,g)$ is irreducible it follows that in the second case $(M,g)$ is
$\Bbb{CP}^2\sharp\overline{\Bbb{CP}}^2$ with D. Page's metric. (see [LeB])$\k$

\medskip
{\bf Remark.  } It is not difficult using the methods from [J-3] to construct Hermitian $\AC$-metrics on $\Bbb{CP}^2$, with two eigenvalues which
 coincide in exactly one point. We shall give the appropriate examples in the last section of the paper.  These metrics are not bi-Hermitian, one of
  the complex structures does not extend to the whole of $\Bbb{CP}^2$, the other one extends to the standard complex structure  on $\Bbb{CP}^2$.
  In fact $\Bbb{CP}^2$ does not admit opposite complex structures.

\medskip
{\bf Proposition 2. } {\it Let us assume that $(M,g)$ is a
bi-Hermitian Gray surface with Hermitian complex structures
$J,\bJ$. If $\xi$ is a Killing vector field on $(M,g)$ such that
$\n_{\xi}J=\n_{\xi}\bJ=0$ then $S_0\xi$ is a Killing vector field,
where $S_0$ is a Killing tensor associated with $\rho$, i.e.
$\rho(X,Y)=g(S_0X,Y)+\frac{\tau}3g(X,Y)$.  }
\medskip
{\it Proof. } Let $S_{\rho}$ be the Ricci endomorphism of $(M,g)$ i.e. $\rho(X,Y)=g(S_{\rho}X,Y)$. Then $S_{\rho}=S_0+\frac{\tau}3$. Since
 $L_{\xi}S_{\rho}=0$ and $L_{\xi}\tau=0$ it is clear that $L_{\xi}S=0$. Both $J,\bJ$ are determined only by $S_{\rho}$ and $g$ thus $L_{\xi}J=L_{\xi}\bJ=0$.
 Thus the result follows from Lemma B.$\k$
\medskip
{\bf Proposition 3. } {\it Let us assume that $(M,g,J)$ is a
compact Hermitian surface with Hermitian Ricci tensor whose group
of (real) holomorphic isometries has a principal orbit of
dimension 3. Then the natural opposite structure $\bJ$ is
Hermitian i.e. complex and orthogonal.}
\medskip
{\it Proof.} Let $\th$ be the Lee form of $(M,g,J)$. Then $|\th|=\frac1{\sqrt2}|\n J|$. If $\zeta$ is a holomorphic Killing vector field then
$\th(\zeta)=0$. It is also clear that $d|\th|^2(\zeta)=0$.
Consequently in an open and dense subset $U$ of $M$ we have $d|\th|^2=f\th$ for some function $f\in C^{\infty}(U)$. The result is now clear in view of
 [J-4], Lemma F.$\k$
\medskip
Our next corollary describes bi-Hermitian Gray surfaces of genus
$0$, i.e. holomorphic $\Bbb{CP}^1$ bundles over $\Bbb{CP}^1$,
which are of cohomogeneity 1 with respect to the group of real
holomorphic isometries.
\medskip
{\bf Corollary.} {\it Let us assume that $(M,g,J)$ is a compact
proper bi-Hermitian Gray surface whose group of (real) holomorphic
isometries has a principal orbit of dimension 3. Then the vector
field $\xi$ coincides with $\eta$ up to a constant factor, the
distribution $\De$ spanned by $\xi,J\xi$ is contained in the
 nullity of both $J,\bJ$ and  $\bJ$ is the natural opposite structure for $J$. The distribution $\De$ coincides with one of eigendistributions of
 the Ricci tensor $S$.}
\medskip
{\it Proof.}  Let us assume that $(M,g)$  is not conformally flat. It means that $|W|\ne 0$. Consequently there exists an open subset $U\in M$ such that
(up to a change of orientation)
$W^-\ne 0$ on $U$. It means that the natural opposite structure for $J$, which is Hermitian in view of Prop.4, coincides in $U$ up to a sign with $\bJ$
 as the only simple eigenvalue of $W^-$. Thus the result of M. Pontecorvo (Prop. 1.3. in [Po]) says that these two structures coincide (up to a sign),
 everywhere where the opposite natural structure to $J$ is defined. Consequently the nullity $\De$ of $J$ coincides with the nullity of $\bJ$ and $\De$
 is one of eigendistributions of the Ricci tensor $\rho$ of $(M,g)$. Since $\th(\xi)=\th(\eta)=0$ it follows  that $\xi=c\eta$ for some $c\in\Bbb{R}-\{0\}$.

If $W=0$ then  $(M,g)$ is conformally equivalent to the product
 $\Bbb{CP}^1\times \Si_g$ where $\Si_g$ is a Riemannian surface of genus $g>0$ and both $\Si_g,\Bbb{CP}^1$ have standard metrics with constant
 opposite sectional curvatures which finishes the proof.$\k$

\medskip
Let us denote by $\bn,
\n,\n^1$ the Levi-Civita connections with respect to the metrics $\bg,g,\bg_1$ respectively.
We have
$$\rho =\brho+2f^{-1}\n df-  f^{-2}(f\de f+3|\n f|^2)\bg,\tag 2.5$$
where $\rho,\brho$ are the Ricci tensors of $(P,g),(P,\bg)$ respectively. The field $\xi=J(\bn f)$ is a holomorphic (with respect to $J$) Killing field
on $(M,g)$ and $(M,\bg)$.  It is easy to see that $\xi= -J\n(\frac1f)$. Analogously the field $\eta=\bJ(\bn^1 h)$ is a holomorphic (with respect to $\bJ$)
 Killing field on $(M,\bg_1)$ and $(M,g)$ and $\eta= -\bJ\n(\frac1h)$.
From Prop.1 it follows that if the scalar curvature $\tau$ of $(M,g)$ is non-constant then both $(M,g,J)$ and $(M,g,\bJ)$ are ruled surfaces.
 Thus $\pi:M\rightarrow \Sigma$ is a holomorphic bundle over a compact Riemann surface $\Si$ with a fiber $\Bbb{CP}^1$. Let us denote by
  $V:=$  ker  $d\pi$ the vertical distribution and by $H=V^{\perp}$ the horizontal distribution of $(M,g)$ induced by the projection
  $\pi:M\rightarrow \Si$ and the metric $g$. Since both structures $J,\bJ$ commute with the Ricci tensor $\rho$ of $(M,g)$ it follows that they
  are determined only by the metric $g$. Consequently every Killing field preserve both structures.  Thus Killing field $\xi$ preserves $\bJ$ and
   $\eta$ preserves $J$, which means that $L_{\xi}\bO=0,L_{\eta}\0=0$.    Now we prove
\medskip
{\bf Proposition 4. }  {\it Let us assume that $(M,g,J,\bJ)$ is a compact bi-Hermitian Gray surface such that $(M,g,J)$ is a ruled surface of genus $g>0$.
Then $\bJ$ is the natural opposite structure for $(M,g,J)$ and the distribution $\De$ spanned by $\xi,J\xi$ is contained in the nullity of  both $J,\bJ$.
 Moreover $\De$ coincides with one of the eigendistributions of the Ricci tensor $S$.}
\medskip
{\it Proof.  }  Let us denote by $S=S_{\rho}$ the Ricci tensor of $\pe$ and by $S_0$ the Killing tensor related with $S_{\rho}$.  Let us recall that ruled
surface different from $\Bbb{CP}^1\times\Bbb{CP}^1$ admits only one ruling (see [B-P-V]). Thus every biholomorphic mapping  $\phi$ must preserve the
 fibers of such a ruled surface, i.e. $\phi(\pi^{-1}(\pi(x)))=\pi^{-1}(y)$ where $y=\pi(\phi(x))$. It follows that the one parameter subgroups of
 holomorphic isometries are $\pi$-related with one parameter subgroups  of biholomorphisms of $\Si$. Thus every holomorphic Killing vector field
 with zeros  $\xi$ on $M$ is $\pi$-related with holomorphic vector field $\xi_0$ with zeros on $\Si$. Consequently if $M$ is of genus $g>0$ then the
 one-parameter subgroups of $\xi,\eta$ both preserve every  fiber of $\pi:M\rightarrow \Si$, i.e.    $\xi,\eta\in \G(V)$.  From Lemma A it follows that
 $\xi=c\eta$ for some constant $c\in\Bbb R-\{0\}$. Consequently $\xi$ belongs to the nullity of both $J,\bJ$, i.e. $\n_{\xi}J=\n_{\xi}\bJ=0$.
 From Prop.2 it follows that $S_0\xi$ is a Killing vector field.
Note that $S_0\xi,\xi,J\xi\in V$ which implies $S_0\xi=\lb_0\xi$. Thus $\lb_0$ is constant.
Since $\xi\in\G(\1)$ and $\1$ is an integrable eigendistribution of a Killing tensor $S_0$ it follows that  $\1$ is totally geodesic (see [J-4] p.7
 Cor.1.4.). Consequently $\bJ$ is the natural opposite structure of $J$ and the distribution $\De$ spanned by $\xi,J\xi$ is the nullity of both
  $J,\bJ$ (see [J-4] Lemma F). In particular $\De$ is $J$ and $\bJ$ invariant, which means that it coincides with one of eigendistributions of $S$.
  On the other hand $\De$ coincides with a vertical distribution $V$ (both have the same section $\xi$ and are $J$ -invariant). Since
  $S_{\rho}\circ J=J\circ S_{\rho}$ and
$S_{\rho}\circ \bJ=\bJ\circ S_{\rho}$ it follows that $V,H$ are eigendistributions of $S$ i.e. $V=\1,H=\2$ where $\lb,\mu$ are eigenvalues of
 $S_{\rho}$ and $\lb=\lb_0+\frac13\tau,\mu=\mu_0+\frac13\tau$. $\k$
\medskip

Since $\eta,\xi$ are Killing fields on $\Bbb{CP}^1$ it follows that $\xi$ has on every fiber exactly two isolated zeros (the north and south poles of
a surface of revolution diffeomorphic to $S^2$.) Let us define $U=\{x\in M:\xi_x\ne 0\}$. Then $U$ is an open and dense subset of $M$.

Our present aim is to prove
\bigskip
{\bf Theorem  1.}  {\it Let us assume that $(M,g,J,\bJ)$ is a compact bi-Hermitian Gray surface of genus $g>0$.  Then $(M,g)$ is locally of
 co-homogeneity 1 with respect to the group of all local isometries of $(M,g)$. The manifold $(U,g)$  is  isometric to the manifold
$ (a,b)\times P_k$
where $(P_k,g_k)$ is a 3-dimensional $\A$-manifold (a circle bundle $p:P_k\rightarrow \Si$ ) over a Riemannian surface $(\Si,g_{can})$ of constant
 sectional curvature  $K\in\{-4,0,4\}$ with a metric
$$g=dt^2+f(t)^2\th^2+h(t)^2p^*g_{can},\tag *$$
where $g_k=\th^2+p^*g_{can}$ and $\th$ is the connection form of
$P_k$ such that $d\th=2\pi k\  p^*\om$, $\om\in H^2(\Si,\Bbb R)$
is an integral, harmonic (hence parallel with respect to
$g_{can}$) 2-form corresponding to the class $1\in H^2(\Si,\Bbb
Z)=\Bbb Z$. The functions $f,h\in C^{\infty}(a,b)$ satisfy the
conditions:

(a) $f(a)=f(b)=0, f'(a)=1, f'(b)=-1, $ ;

(b)   $h(a)\ne0\ne h(b),h'(a)=h'(b)=0, $ .}
\medskip
{\it Proof.}  The best way to prove this theorem is to use the recent results contained in [A-C-G]. For a while we shall use a notation from [A-C-G].
 Note that Proposition 5 yields that $J\xi=\bar J\xi$. Consequently for both metrics $g,\bg$ the natural opposite structure for $\bar J$ coincides with
  $J$. It implies that the K\"ahler surface $(M,\bg,J)$ is of Calabi type (see [A-C-G]). Note that the conformal factor $f$ to the K\"ahler metric
  is the square of an affine function of  the momentum map $z$ of $\xi$ with respect to $\bar\om(X,Y)=\bg(JX,Y)$. The scalar curvature of both these
   metrics is a function of the momentum map $z$ alone (see Prop.5 below and note that $\n\lb=\frac13\n\tau,\n\mu=-\frac16\n\tau$). One can also easily
   check using [J-4] that conformal scalar curvature $\kappa$ and functions $\a,\beta$  also depend only on  $z$.
It follows from [A-C-G], Lemma 10, that the scalar curvature $s_{\Si}$ is constant. Thus it follows from the methods of Lebrun
(see Prop.13 in [A-C-G]) that
both  metrics $g,\bg$ are local cohomogeneity 1. Consequently on  the open, dense subset, where $\xi\ne 0$, the metric $g$ is of the form $(*)$.
The boundary conditions are  the conditions (a),(b) in view of [B],[M-S].$\k$
\medskip
{\it Remark } Note that $\om$ depends only on the complex structure $J$  of a Riemannian surface $\Si$. The complex structure $J$ determines a conformal
class of a Riemannian metric $[g]$ such that $g(JX,JY)=g(X,Y)$. The matric $g_{can}$ is the metric in this class of constant sectional curvature.
Consequently if $\Si$ is a Riemannian surface of genus $g$ then every complex structure $J$ on $\Si$ determines a unique form $\om$ and consequently a
 family of $S^1$-principle bundles $P_{k,\Si}$.
Note also  that it is not true in general that every local bi-Hermitian $\AC$-metric of non-constant scalar curvature is local cohomogeneity one metric.
The counterexample gives the
(non-compact) Einstein-Hermitian self-dual space $(M,g)$ of co - homogeneity 2 constructed by Apostolov and Gauduchon in [A-G-2],Th.2. The related
 K\"ahler metric $(M,\bg)$ is weakly self-dual of co-homogeneity grater than 1. In fact every Killing vector field with respect to $(M,\bg)$ is also
 a Killing vector field for $(M,g)$ ( the conformal factor to an Einstein metric is the square of the scalar curvature $\bt$ of $(M,\bg)$ - see [D-1],
 Prop.4).
\bigskip
We shall end this section with characterization of the eigenvalues
of a bi-Her\-mitian Gray surface of genus $g>0$. We show that Lie
forms $\th,\bth$ and the difference $\lb-\mu$ of eigenvalues of
the Ricci tensor of a  Hermitian Gray surface $(M,g,J,\bJ)$
 depend only on the length of tensor fields $\n J,\n \bJ$.
\medskip
{\bf Proposition  5.}  {\it Let us assume that $(M,g,J,\bJ)$ is a compact bi-Hermitian Gray surface such that $\bJ$ is the natural opposite
 Hermitian structure of $J$.  Let $\th,\bth$ be the Lee forms of $(M,g,J)$ and $(M,g,\bJ)$ respectively. Then
$$\gather\th=2d\ln\frac{\gamma}{|1-\e \gamma|},\tag 2.6a\\
\bth=-2d\ln|1-\e \gamma|,\tag 2.6b\\
\lb-\mu=C\frac \gamma{(1-\e \gamma)^2},\tag 2.6c\endgather$$
where $\e\in\{-1,1\},C\in \Bbb{R}-\{0\},\gamma=\frac{\be}{\a},\a=|\n J|,\be=|\n\bJ|$.}

\medskip
{\it Proof. } From [J-4] it follows that
$$-d\ln\a-\frac12\th=-d\ln\be-\frac12\bth.\tag 2.7$$
Consequently we obtain:
$$\gather
\th-\bth=2d\ln\frac{\be}{\a},\tag 2.8a\\
\th+\bth=2d\ln|\lb-\mu|.\tag 2.8b\endgather$$
Thus
$$
\th=d\ln\frac{\beta|\lb-\mu|}{\a},
\bth=d\ln\frac{\a|\lb-\mu|}{\be}.\tag 2.9$$
On the other hand (see [J-4])
$$\beta\th=\e\a\bth,$$
for a certain $\e\in\{-1,1\}$.
Consequently
$$\n\ln|\lb-\mu|=-(\frac{\e+\gamma}{\gamma-\e})\frac{d\gamma}{\gamma}.$$
It implies
$$\lb-\mu=C\frac \gamma{(1-\e \gamma)^2},$$
for a certain $C\in\Bbb{R}-\{0\}$. Now it is clear that all formulas $2.6$ hold true.$\k$

\bigskip

{\bf 3. Bi-Hermitian Gray surfaces with genus $g\ge 1$ .}  In this
section we shall construct bi-Hermitian metrics $g$ on ruled
surfaces $(M_{k,g},g)$
 of genus $g$.  Then, according to Th.1., $(M_{k,g},g)$ is locally of co-homogeneity 1 with respect to the group of all local isometries of
 $(M_{k,g},g)$ and an open, dense submanifold $(U_{k,g},g)\subset (M_{k,g},g)$  is  isometric to the manifold
$ (a,b)\times P_k$
where $(P_k,g_k)$ is a 3-dimensional $\A$-manifold (a circle bundle $p:P_k\rightarrow \Si_g$ ) over a Riemannian surface $(\Si_g,g_{can})$ of
constant sectional curvature  $K\in\{-4,0,4\}$ with a metric
$$g_{f,g}=dt^2+f(t)^2\th^2+g(t)^2p^*g_{can},\tag 3.1$$
where $g_k=\th^2+p^*g_{can}$ and $\th$ is the connection form of $P_k$ such that $p^*d\th=2\pi k\  \om$, $\om\in H^2(\Si_g,\Bbb R)$ is an integral
form, parallel with respect to $g_{can}$, corresponding to the class $1\in H^2(\Si_g,\Bbb Z)=\Bbb Z$. It follows that $P_1=G\slash \G$, where $\G$
is a lattice in $G=\widetilde{SL(2,\Bbb{R})},G=H $ or $G=SU(2),\G=\{e\}$ and $P_k=\Bbb{Z}_k\backslash G\slash\Gamma$. Let $\th^{\sharp}$
be a vector field dual to $\th$ with respect to $g_P$. Let us consider a local orthonormal frame $\{X,Y\}$ on $(\Si_g,g_{can})$ and let $X^h,Y^h$
be horizontal lifts of $X,Y$ with respect to $p:M_{k,g}\rightarrow \Si_g$ (i.e. $dt(X^h)=\th(X^h)=0$ and $p(X^h)=X$) and let
$H=\frac{\partial}{\partial t}$.
Let us define two almost Hermitian structures $J,\bJ$ on $M$ as follows
$$JH=\frac1f\th^{\sharp},JX^h=Y^h,\ \bJ H=-\frac1f\th^{\sharp},\bJ X^h=Y^h.$$
\medskip
{\bf Proposition 6. } {\it Let $\De$ be a distribution spanned by the fields $\{\th^{\sharp},H\}$.
Then $\De$ is a totally geodesic foliation with respect to the metric $g_{f,g}$. Both structures $J,\bJ$ are Hermitian and $\De$ is contained in
 the nullity of $J$ and $\bJ$.
The distribution $\DE$ is umbilical with the mean curvature normal $\xi=-\n\ln g$.  Let $\lb, \mu$ be eigenvalues of the Ricci tensor $S$ of
$g_{f,g}$ corresponding to eigendistributions $\De,\De^{\perp}$ respectively.  Then the following conditions are equivalent:

(a) There exists $E\in \Bbb R$ such that $\lb-\mu =Eg^2$,

(b) There exist $C,D\in \Bbb R$ such that $\mu =Cg^2+D$,

(c) $\lb-2\mu$ is constant,

(d) $(U_{k,g},g_{f,g})$ is a bi-Hermitian Gray surface.}
\medskip
{\it Proof.  } The first part our Proposition is a consequence of [J-2]. Note that $\n\lb=H\lb H,\n\mu=H\mu H$. Consequently
$tr_g\n S=\frac12\n\tau=(H\lb+H\mu)H$. On the other hand one can easily check that $tr_g\n S=2(\mu-\lb)\xi+H\lb H$. Thus

$$\frac{\n\mu}{2(\lb-\mu)}=\n\ln g.$$

Now we prove that (a) $\Rightarrow$ (b). If (a) holds then $\n \mu=2Eg^2\frac{\n g}g=E\n g^2$. Thus $\n(\mu-Eg^2)=0$ which implies (b).

(b)$\Rightarrow$ (a). We have
$$-\frac{\n g}g=\frac{\n\mu}{2(\mu-\lb)}=\frac{Cg\n g}{\mu-\lb},$$
and consequently $\n g(\frac{Cg^2+\mu-\lb}{g(\mu-\lb)})=0$ which is equivalent to (b).

(a)$\Rightarrow $(c). We have $\lb-\mu=Eg^2$ and consequently  $\n \mu=2Eg\n g=E\n g^2$. Thus $\n \lb=\n (\mu+Eg^2)=2E\n g^2$ and
$\n\lb-2\n\mu=0$ which gives (c).

(c)$\Rightarrow$ (a). If $\n\lb=2\n\mu$ then $\n\lb=4(\lb-\mu)\frac{\n g}g$. Consequently $\n\lb-\n\mu=2(\lb-\mu)\frac{\n g}g$ and
$\n\ln|\lb-\mu|=2\frac{\n g}g=2\n\ln g$, which means that $\n\ln|\lb-\mu|g^{-2}=0$.  It follows that $\ln\frac{|\lb-\mu|}{g^2}=C$ for
some $C\in\Bbb R$, which is equivalent to  (a).

(d)$\Leftrightarrow$(c).  This equivalence follows from [J-3].
$\k$
\medskip
{\bf Theorem  2 } {\it On any  ruled  surface $M_{k,g}$ of genus $g>0$ with $k>0$  there exist a  one-parameter family of   Hermitian $\AC$-metrics
 $\{g_{x}:x\in(0,1)\}$ which contains all bi-Hermitian Gray metrics on $M_{k,g}$.}

\medskip
{\it Proof.}
 Note that for the first Chern class $c_1(\Si_g)\in H^2(\Si_g,\Bbb Z)$ of the complex curve $\Si_g$ we have the relation $c_1(\Si_g)=\chi\a$,
 where $\a\in H^2(\Si_g,\Bbb Z)$ is an indivisible integral class and $\chi=2-2g$ is the Euler characteristic of $\Si_g$.
 Let us write $s=\frac{2k}{|\chi|}$ if $g\ne 1$ and $s=k$ if $g=1$. Then it is easy to show that the manifold $(M_{k,g},g)$ with the metric $g$
 given by $(*)$ has the Ricci tensor with the following eigenvalues :
$$\gather
\lb_0=-2\frac{g''}g-\frac{f''}f,\tag 3.1a\\
\lb_1= -\frac{f''}f-2\frac{f'g'}{fg}+2s^2\frac{f^2}{g^4},\tag 3.1b\\
\lb_2=
-\frac{g''}g-\frac{f'g'}{fg}-(\frac{g'}g)^2-2s^2\frac{f^2}{g^4}+\frac{K}{g^2},\tag
3.1c
\endgather$$
where $\lb_0,\lb_1,$ correspond to eigenfields $T=\frac d{dt},\th^{\sharp}$ and $\lb_2$ corresponds to a two-dimensional eigendistribution orthogonal to
$T$ and $\th^{\sharp}$.  If $(M,g)\in\AC$ is a bi-Hermitian Gray surface then $\lb_0=\lb_1=\lb$ and, if we  denote $\mu=\lb_2$, Prop.6 and [J-3] imply
an equation $$\mu=Dg^2-C\tag3.2$$
 for some $D,C\in \Bbb R$. Since $\lb_0=\lb_1$ we get
$$f=\pm\frac{gg'}{\sqrt{s^2+Ag^2}}.\tag3.3$$
 Using a homothety of the metric we can assume that  $A\in\{-1,0,1\}$. In the case $A= 0$ we get a weakly-self-dual K\"ahler metric and these metrics
 on compact complex surfaces are classified (see [A-C-G]). So we restrict our considerations to the case $A\in\{-1,1\}$.
 Now we introduce a function $h$ such that $h^2=s^2+Ag^2$. Note that im$  h\subset (-s,s)$ if $A=-1$ and im$  h\subset (s,\infty)$.
 if $A=1$.
Then $g=\sqrt{|s^2-h^2|}$. Let us introduce a function  $z$, such that $h'=\sqrt{z(h)}$.  Note that
$$ f=h'\text{   and    } f'=\frac12z'(h).\tag 3.4$$
It follows that equation (3.2) is equivalent to
$$z'(h)-z(h)\frac{s^2+h^2}{h(s^2-h^2)}=\frac{4\e}h+\frac{D(s^2-h^2)^2}h-\frac{C(s^2-h^2)}h,\tag3.5$$
where $\e =\text{sgn} KA\in\{-1,0,1\}$.  It follows that
\define \h{(\frac hs)}
$$\gather  z(h)=(1-\h^2)^{-1}(-4\e\h^{2}-\frac {Ds^4}5\h^6+(Ds^4-\frac{Cs^2}3)\h^4+\tag 3.7\\+(2Cs^2-3Ds^4)\h^2-4\e+Cs^2-Ds^4+\frac Es \frac hs).
\endgather  $$

Let us denote again $C=Cs^2,D=Ds^4,E=\frac Es$ and let
$$z_0(t)=(1-t^2)^{-1}(-4\e(1+ t^{2})+D(-\frac15t^6+t^4-3t^2-1)+C(-\frac13
t^4+2t^2+1)+Et)
\tag 3.8$$
Write
$$P(t)=(-4\e t^{2}-\frac{D}5t^6+(D-\frac{C}3)t^4+(2C-3D)t^2+Et-4\e+C-D).\tag 3.9$$
Then $z_0(t)=\frac{P(t)}{1-t^2} $. Note that  $z(h)=z_0\h$ and $z'(h)=\frac1sz_0'\h$. We are looking  for real numbers $x>y\in \Bbb R$ such that
$$\gather
z_0(x)=0,  z'_0(x)=-2s,\tag 3.10a \\
z_0(y)=0,  z'_0(y)=2s,\tag 3.10b \endgather$$
and $z(t)>0$ for $t\in(y,x)$.
  Note that equations (3.10a) are equivalent to
$$\gather
-4\e x^{2}-\frac{D}5x^6+(D-\frac{C}3)x^4+(2C-3D)x^2-4\e+C-D+Ex=0\tag 3.11a\\
-8\e x-\frac{6D}5x^5+4(D-\frac{C}3)x^3+2(2C-3D)x+E=-2s(1-x^2).\tag3.11b\endgather$$
Equations $(3.11)$ yield
$$\gather D=\frac{5(-3E-6s-24\e x+3Ex^2-12sx^2-8\e x^3+2s x^4)}{2(-1+x)x(1+x)(15+10x^2-x^4)},\tag 3.12a\\
C=\frac{3(5E+10s+80\e x+30 s x^2-10E x^2+5Ex^4-10sx^4-16\e x^5+2s x^6)}{2(-1+x)x(1+x)(-15-10x^2+x^4)}\tag 3.12b
\endgather$$
Solving in a similar way equations $(3.10b)$ one can see that there exists a function $z_0$ satisfying the equations (3.10) if
$$\gather
(x+y)(-4\e(-5 x+x^3+5 y+2x^2y-2xy^2-y^3)\tag 3.13\\+s(5+2x^3y+2xy^3+3y^2+3x^2+x^2y^2-16xy))=0,\endgather$$
where  $x>y$, $x,y\in (-1,1)$ in the case $A=-1$ and $x,y\in (1,\infty)$ in the case $A=1$.
Using standard methods one can check that in the case of the genus $g\ge 1$ (i.e. if $K=-4$ or $K=0$)
the only solutions of (3.13) giving a positive function $z$ are $x=-y\in(0,1)$. In the case
$g=0,K=4$ apart from the solutions with $x=-y$ (see [J-3]) there are two additional families of solutions with $\e=1$ and $\e=-1$ on the first
Hirzebruch surface $F_1$ and one additional family with $\e=-1$   on the second Hirzebruch surface $F_2$.

It follows that if $g\ge 1$ then $x=-y,E=0$ and $\e=1$ or $\e=0$.
Consequently
$$\gather
P(t)\tag 3.14\\=\frac1{x(15-5x^2-11x^4+x^6)}((t^2-x^2)(s(-15+10x^2-3x^4+t^2(10+12x^2-6x^4)\\
+t^4(-3-6x^2+x^4))+4\e x(x^2(-5+x^2)-t^4(3+x^2)+t^2(5+2x^2+x^4)))).\endgather$$
Thus $$P(0)=\frac{-4x^4(x^2-5)+sx(15-10x^2+3x^4)}{15-5x^2-11x^4+x^6}$$ and
$P(t)>0$ if $t\in(0,x)$ for all $x\in(0,1)$.
  Now the function $z_0(t)=\frac1{(1-t^2)}P(t)$ is positive on $(-x,x),x\in(0,1)$. If $x\in(0,1)$ then there exists a solution
  $h :(-a,a)\rightarrow (-sx,sx)$, where
$$a=\lim_{t\rightarrow sx^-}\int^t_0\frac {dh}{\sqrt{z_0\h}},$$
of an equation
$$h'=\sqrt{z_0\h},$$
such that $h(-a)=-sx,h(a)=sx,h'(-a)=h'(a)=0,h''(-a)=1,h''(a)=-1$. It follows that functions $f=h',g=\sqrt{s^2-h^2}$ are smooth on $(-a,a)$
 and satisfy the boundary conditions described in Th.2. Consequently the metric
$$g_x=dt^2+f(t)^2\th^2+g(t)^2p^*g_{can}, $$
on the manifold $(-a,a)\times P_k$ extends to the smooth metric on the compact ruled surface $M=P_k\times_{S^1}S^2$ which is a $2-$sphere bundle
 over Riemannian surface $\Si_g$. Note that $g(-a)=g(a)=s\sqrt{1-x^2}$ and that our construction is valid for all $x \in(0,1)$. $\k$
\medskip
{\bf Theorem  3. } {\it There are no irreduciblr bi-Hermitian Gray-metrics  on trivial ruled surfaces $M_{0,g}=\Bbb{CP}^1\times \Sigma_g$ with $g>0$.}

\medskip
{\it Proof.}
Now we consider Gray metrics on the product $\Bbb{CP}^1\times \Sigma_g$.  Then, $K=-4,s=0$ and
we can take $g=h,f=h'$. Consequently
$z(h)=-4+\frac{D}5h^4+\frac C3h^2+\frac Eh$. For simplicity let us write $D=\frac D5,C=\frac C3$.
Then $z(h)=-4+Dh^4+Ch^2+\frac Eh$. We are looking for solutions of an equation $h'=\sqrt{z(h)}$ satisfying initial conditions $h(a)=x,h(b)=y$
 where  $$z(x)=0,z'(x)=2,z(y)=0,z'(y)=-2\tag 3.15$$  for some unknown $x,y\in\Bbb R$ such that $0<x<y$.
Equation (3.15) yields
$$\gather
D=-\frac4{x^2y^2}+\frac{E(x^3-y^3)}{x^3y^3(x^2-y^2)},C=E\frac{y^5-x^5}{x^3y^3(x^2-y^2)}+4\frac{x^2+y^2}{x^2y^2},\tag 3.16a\\
D=E\frac{y^3-x^3}{4x^3y^3(x^2-y^2)}+\frac1{2xy(x-y)},C=E\frac{x^5-y^5}{2x^3y^3(x^2-y^2)}-\frac{x^3+y^3}{xy(x^2-y^2)}.\tag 3.16b\endgather$$
It follows that equations (3.16) have a solution if and only if
$$E=\frac25\frac{xy(y+x)}{x^3-y^3}(8(x-y)+xy)=\frac23\frac{xy}{x^5-y^5}(4(x^4-y^4)+xy(x^3+y^3)).$$
If $y=\a x$ where $\a \ne 1$ then we obtain $$
x=\frac{-4(\a-1)(\a^2+3\a+1)}{\a(\a^2+\a+2)}.$$
Since $0<x<y$ and $\a>1$ we get a contradiction.
 Consequently there are no irreducible Gray metrics on the trivial ruled surface $M=\Bbb{CP}^1\times \Sigma_g$,  where $g>1$.

We shall finish by investigating Gray metrics on the surface $\Bbb {CP}^1\times T^2$. Now
$z(h)=\frac{D}5h^4+\frac C3h^2+\frac Eh$. For simplicity let us write $D=\frac D5,C=\frac C3$.
Then $z(h)=Dh^4+Ch^2+\frac Eh$. We are looking for solutions  satisfying initial conditions (3.15).
It follows analogously as above that equations (3.15) have a solution if and only if
$$E=\frac25\frac{x^2y^2(y+x)}{x^3-y^3}=\frac23\frac{x^2y^2(x^3+y^3)}{x^5-y^5}.\tag 3.17$$
Consequently, if $y=\a x$ then  we obtain $(\a+1)(\a-1)^3(2\a^2+\a+2)=0,$ where $\a>1$. It follows that there are no irreducible Gray metrics on
 the trivial ruled surface $M=\Bbb{CP}^1\times T^2$.$\k$

\bigskip
{\bf  4.  Hermitian Gray structures on $\Bbb{CP}^2$. }  In this section we give examples of $\AC$-4-manifolds $(M,g)$ whose Ricci tensor $\rho$
 has two eigenvalues $\lb,\mu$ such that only one from two natural complex structures defined by the Killing tensor $\rho-\frac13\tau g$ on the
  subset $U=\{x\in M:\lb(x)\ne\mu(x)\}$ extends to the complex structure on the whole of the manifold $M$. Let us denote by $J$ the standard
  complex structure of a projective space $\Bbb{CP}^2$.
As in [M-S], [J-3] by $L(k,1)$ (where $k\in\Bbb{N}$) we shall denote the Lens spaces.  The manifolds $L(k,1)$ are the circle bundles over $\Bbb{CP}^1$.
 Then $\Bbb{CP}^2$  is the space of cohomogeneity 1 under an action of $U(2)$ with principal orbit $P=L(k,1)$ (with $k=1$) and two special orbits:
 $\Bbb{CP}^1$ and a point (i.e.   $\Bbb{CP}^2=[\Bbb{CP}^1|S^3|*]$). Let us denote by $\eta$ the only real eigenvalue of the polynomial
 $S(x)=x^3+5x^2+75x+59$. Then  $\eta=-\frac53-\frac{100\root {3}\of{4}}{3(383+129\sqrt{129})}+\root{3}\of{\frac13(2(383+129\sqrt{129}))}=-0.8245....$.
Now our aim is to prove the following theorem:
\medskip
{\bf Theorem  4. } {\it On the  surface $\Bbb {CP}^2$ there exist two  one-parameter families of   Hermitian $\AC$-metrics
 $\{g_{x}:x\in T=(\eta,1)\cup (1,\infty)\}$. The Ricci tensor $\rho=\rho_{x}$ of $(\Bbb{CP}^2,g_{x})$ is Hermitian with respect to
 the standard complex structure of $\Bbb{CP}^2$ and has two eigenvalues, which coincide in an exactly one point. Every
co-homogeneity one $\AC$-metric  with $J$-invariant Ricci tensor on  $\Bbb{CP}^2$ is homothetic to one of $g_x$ where  $x\in T$.}
\bigskip

{\it Proof. } We shall retain the notation of section 3. Note that the genus $g=0$ and all
formulas (3.1)-(3.7) remain valid for this case, however the boundary conditions will be different.
We shall find the conditions on $f,g$ to extend  the  metric $(3.1)$ on the whole of $\Bbb{CP}^2$. The metric (3.1) extends to the metric on
$[\Bbb{CP}^1|S^3|*]$  if and only if
the boundary conditions are as follows:
  $$\gather f(b)=g(b)=0,f'(b)=g'(b)=-\e, \tag 4.1a\\
g(a)\ne 0,g'(a)=0,f(a)=0,f'(a)=1.\tag 4.1b\endgather $$
At first we shall consider $(4.1a)$. We get $h(b)=1$ and consequently  $z_0(1)=0,z_0'(1)=-2k\e$ where $z_0(t)=\frac{P(t)}{1-t^2}$
and $P(t)=(-4\e(1+ t^{2})+D(-\frac15t^6+t^4-3t^2-1)+C(-\frac13t^4+t^2+1)+Et)$. Equations (4.1a) imply $E=-\frac8{15}(5C-6D+15\e)$.
Then $z'_0(1)=-2\e$ and consequently $k=1,z=z_0$.
Now we consider (4.1b).
We have to find $x\in(-1,1)$ if $\e=1$ or $x\in(1,\infty)$ if $\e=-1$  such that $g(a)=x, z(x)=0,z'(x)=2\e$.
Then $f(a)=0,f'(a)=1,g(a)=\sqrt{|1-x^2|}\ne 0, g'(a)=0$. These equations are equivalent to
$$D=\frac{-5\e}{4+x-4x^2-x^3},\ C=\frac{3\e(7+4x-x^2)}{(-1+x)(1+x)(4+x)}.$$
Consequently
$$\gather
z_x(t) = \\
\frac{\e(t-1)(t-x)(t^3+t^2(2+x)+t(5+6x)+8+13x+4x^2))}{(1+t)(x-1)(1+x)(4+x)}.\endgather$$
Using elementary calculations one can prove that the polynomial $Q_x(t)=(t^3+t^2(2+x)+t(5+6x)+8+13x+4x^2))$ is positive for $t\in(x,1)$
 where $x<1$ if and only if $x\in (\eta,1)$ where
  $\eta=-\frac53-\frac{100\root {3}\of{4}}{3(383+129\sqrt{129})}+\root{3}\of{\frac13(2(383+129\sqrt{129}))}=-0.8245....$
  is the only real eigenvalue of the polynomial $S(x)=x^3+5x^2+75x+59$.
Now it is clear that if  $x\in(\eta,1)\cup(1,\infty)$ then the function  $z_x(t)$ has exactly two roots $x,1$ in one of the  intervals
$[x,1]$ and $[1,x]$ respectively and in both considered cases $z_x(t)>0$ if $t\in(x,1)$ or $t\in(1,x)$ respectively. If h is a solution of an equation
$h'=z_x(h)$ satisfying the boundary conditions $h(a)=x,h(b)=1$ then it is easy to verify that functions $f=h'$  and $g=\sqrt{|1-h^2|}$
are positive in $(a,b)$ and satisfy equations $(3.2)$ and boundary conditions $(4.1)$. Consequently our metric defined on $(a,b)\times S^3$
extends to the Gray metric on $\Bbb{CP}^2$. From the construction it is also clear that the Ricci tensor of $(\Bbb{CP}^2,g_x)$ is invariant
with respect to the standard complex structure $J$ of $\Bbb{CP}^2$ and that the opposite complex structure id defined everywhere except
the point corresponding to the degenerate orbit $*$.

\bigskip
\centerline{\bf References.}
\par
\medskip
\cite{A-C-G} V. Apostolov, D. Calderbank and P. Gauduchon {\it The geometry of weakly selfdual K\"ahler surfaces} Compositio Math. 135: 279-322,2003.
\par
\medskip
\cite{A-G-1} V. Apostolov and P. Gauduchon {\it Self-dual Einstein Hermitian four-manifolds} Ann. Scoula Norm. Sup. Pisa Cl. Sci (5), vol.I(2002),203-243.
\par
\medskip
\cite{A-G-G} V. Apostolov, P. Gauduchon and P. Grantcharov {\it Bi-Hermitian structures in complex surfaces} Proc. London Math. Soc. (3) 79 (1999), 414-428.

\par
\medskip
\cite{A-G-2} V. Apostolov and P. Gauduchon {\it The Riemannian Goldberg-Sachs theorem} Int. J. Math. 8,(1997),421-439

\par
\medskip
\cite{B} B\'erard-Bergery, L.{\it Sur de nouvelles vari\'et\'es riemanniennes d'Einstein}, Institut \'Elie Cartan,6, Univ. Nancy, (1982), 1-60.
\par
\medskip
\cite{Be}  A. Besse, {\it Einstein manifolds}, Springer Verlag, 1987.
\par
\medskip
\cite{Bo}  C.P.Boyer, {\it Conformal duality and compact Complex Surfaces}  Math. Ann. 274, 517-526 (1986).

\par
\medskip
\cite{D-1} A. Derdzi\'nski ,{\it Self-dual K\"ahler manifolds and Einstein manifolds of dimension four}, Compositio Math. 49, (1983), 405-433.
\par
\medskip
\cite{D-2} A. Derdzi\'nski ,{\it Exemples de m\'etriques de K\"ahler et d'Einstein auto-duales sur le plan complexe}, G\'eom\'etrie riemannienne en dimension 4 ( S\'eminaire Arthur Besse 1978-1979) Cedic-Fernand Nathan, Paris (1981), pp. 334-346.
\par
\medskip
[G] A. Gray {\it Einstein-like manifolds which are not Einstein} Geom. Dedicata {\bf 7} (1978) 259-280.
\par
\medskip
\cite{J-1} W. Jelonek,  {\it Compact K\"ahler surfaces with harmonic anti-self-dual Weyl tensor}, Diff. Geom. and Appl. 16, (2002), 267-276.
\par
\medskip
\cite{J-2} W. Jelonek,  {\it Extremal K\"ahler $\AC$-surfaces}, Bull. Belg. Math. Soc. 9 (2002), 561-571.
\par
\medskip
\cite{J-3} W. Jelonek,  {\it Bi-hermitian Gray surfaces},  Pacific
J. Math. 222  no. 1,  57-68 , 2005
\par
\medskip
[J-4] W. Jelonek {\it Einstein Hermitian and anti-Hermitian surfaces  }   Ann. Polon. Math. 81.1 (2003) ,7-24.
\par
\medskip
\cite{K} T. Koda, {\it A remark on the manifold $\Bbb {CP}^2\sharp\overline{\Bbb {CP}}^2$ with B\'erard-Bergery's metric},  Ann. of Global Analysis
and Geom. {\bf 11}  (1993), 323-329.
\par
\medskip
\cite{LeB} C. LeBrun {\it  Einstein metrics on complex surfaces}, in Geometry and Physics (Aarhus 1995) eds. J. Andersen, J. Dupont, H. Pedersen and
A. Swann, Lect. notes in Pure Appl. Math., Marcel Dekker, 1996.

\par
\medskip
\cite{P} D. Page, {\it A compact rotating gravitational
instanton}, Phys. Lett.79 {\bf B} (1978), 235-238.
\par
\medskip
[M-S] B.Madsen, H. Pedersen, Y. Poon, A. Swann {\it Compact Einstein-Weyl manifolds with large symmetry group.} Duke Math. J. {\bf 88} (1997), 407-434.

\par
\medskip
\cite{Po} M. Pontecorvo , {\it Complex structures on Riemannian
four-manifolds} Math. Ann. 309 (1997), no. 1, 159-177.
\par
\medskip
\cite{S} P. Sentenac ,{\it Construction d'une m\'etriques
d'Einstein  sur la somme de deux projectifs complexes de dimension
2}, G\'eom\'etrie riemannienne en dimension 4 ( S\'eminaire Arthur
Besse 1978-1979) Cedic-Fernand Nathan, Paris (1981), pp. 292-307.
\par
\medskip
\cite{V} I. Vaisman {\it On locally and globally conformal K\"ahler manifolds } Trans.Amer. Math. Soc {\bf 262} (1980), 533-542.
\noindent

Institute of Mathematics

Technical University of Cracow

 Warszawska 24

31-155 Krak/ow,POLAND.

E-mail address: wjelon\@usk.pk.edu.pl

\enddocument